\documentclass[a4paper,12pt]{article}
\usepackage[top=3cm,bottom=3cm,left=2.5cm,right=2.5cm]{geometry}
\usepackage{t1enc}
\usepackage[latin1]{inputenc}
\usepackage[english]{babel}
\usepackage{amsmath,amsthm}
\usepackage{amsfonts}
\usepackage{latexsym}
\usepackage[dvips]{graphicx}
\usepackage{graphicx}
\usepackage{color}
\usepackage{amsmath}
\usepackage{amssymb}
\usepackage{amsbsy}
\usepackage{amsfonts}
\usepackage{amsfonts,mathrsfs}

\begin{document}

\newtheorem{prop}{Proposition}[section]
\renewcommand{\theprop}{\arabic{section}.\arabic{prop}}
\newtheorem{lem}[prop]{Lemma}
\renewcommand{\thelem}{\arabic{section}.\arabic{lemma}}
\newtheorem{thm}[prop]{Theorem}
\renewcommand{\thethm}{\arabic{section}.\arabic{thm}}
\newtheorem{cor}[prop]{Corollary}
\renewcommand{\thecor}{\arabic{section}.\arabic{Cor}}
\newtheorem{conj}[prop]{Conjecture}
\renewcommand{\theconj}{\arabic{section}.\arabic{conj}}
\renewcommand{\theequation}{\arabic{section}.\arabic{equation}}
\renewcommand{\thefigure}{\arabic{figure}}
\newtheorem{exa}{prop}
\renewcommand{\theexa}{\arabic{section}.\arabic{ex}}
\newtheorem{rem}[prop]{Remark}
\renewcommand{\therem}{\arabic{section}.\arabic{rem}}
\renewcommand{\theequation}{
\arabic{equation}}

\newtheorem{defi}[prop]{Definition}
\renewcommand{\thedefi}{\arabic{section}.\arabic{defi}}
\renewcommand{\thefigure}{\arabic{figure}}
\def\zm{\noindent{\bf Proof.\ }}
\def\ezm{\vspace*{6mm}\framebox{} }

\addtolength{\baselineskip}{+0.6mm}

\title{Leap eccentric connectivity index of some graph operations with subdivided edges
\thanks{This work is supported by the Department of Education of Hunan Province(19A318).}}
\author{Ling Song, Zikai Tang\thanks{Corresponding author: zikaitang@163.com}
 \\
{\small School of Mathematics and Statistics, Hunan Normal University,}\\ {\small Changsha, Hunan 410081, P. R. China.}\\
}

\date{}
\maketitle
{\bf Abstract}: The leap eccentric connectivity index of $G$ is defined as $$L\xi^{C}(G)=\sum_{v\in V(G)}d_{2}(v|G)e(v|G)$$ where $d_{2}(v|G) $ be the second degree of the vertex $v$ and $e(v|G)$ be  the eccentricity of the vertex $v$ in $G$. In this paper, we first give a counterexample for that if $G$ be a graph and $S(G)$ be its  the subdivision graph, then each vertex $v\in V(G)$, $e(v|S(G))=2e(v|G)$ by Yarahmadi in \cite{yar14} in Theorem 3.1.  And we describe the upper and lower bounds of the leap eccentric connectivity index of four graphs based on subdivision edges, and then give the expressions of the leap eccentric connectivity index of join graph based on subdivision, finally, give the bounds of the leap eccentric connectivity index of four variants of the corona graph.


\baselineskip=0.30in

\section{Introduction}
In this paper, we only consider simple, undirected, finite graphs. Let $G$ denote a graph with $n$ vertices and $m$ edge sets. Denoted by $d_{G}(u,v)$ the shortest path of length connecting $u$ and $v$ in $G$, for vertices $u$, $v\in V(G)$. For a vertex $v$ and a positive integer $k$, we let $N_{k}(v|G)$ denote the open $k-neighborhood$ of vertex $v$ in $G$ and defined as $N_{k}(v|G)=\{u\in V(G)|d_{G}(u,v)=k\}$. Let $d_{k}(v|G)$ denote the $k$ degree of the vertex in $G$, expressed as the number of vertices in the open $k-neighborhood$ of vertex $v$ in $G$, that is, $d_{k}(v|G)=|N_{k}(v|G)|$. We can see that for any vertex $v$ in $G$ there are $d_{1}(v|G)=|N_{1}(v|G)|$ and $d_{2}(v|G)=|N_{2}(v|G)|$. The eccentricity is defined as $e(v|G)$, for a vertex $v$ in $G$, which represents the maximum distance from vertex $v$ to other vertices in the graph, that is, $e(v|G)=max\{d_{G}(u,v)|u\in V(G)\}$. For any vertex in the graph, we define the maximum eccentricity value as the diameter $diam(G)$. We let $V^{\alpha}_{e}(G)\subseteq V(G) $ denote the set of vertices in $G$ where the eccentricity is equal to $ \alpha$, where $ \alpha=1,2,...,diam(G)$, obviously $V_{e}^{1}(G)$ represents the set of vertices in $G$ that have a eccentricity of 1 to other vertices, the degree of these vertices is $n-1$, we call these vertices full vertices. Let $H\subseteq V(G)$ denote any subset of vertices of $G$, then the induced subgraph $\langle H\rangle$ of $G$ is the graph that the vertex set is  $H$, and the edge set is the edge in graph $G$ with the vertex in $H$ as the endpoint. If there are no graphs isomorphic to graph $F$ in all induced subgraphs of graph $G$, we call graph $G$ the $F-free$ graph. For other terms and symbols that are not defined here, please refer to the reference \cite{bod08}.

Structure descriptors based on molecular graphs are often called topological indices and have very important meanings. In 1972, Gutman and Trinajestic \cite{gut72} introduced  the first, and elaborated it in \cite{gut75}. The definition is as follows:
$$ M_{1}(G)=\sum_{v\in V(G)}d_{1}^{2}(v|G). $$

Recently, some scholars replaced the vertex degree in the first Zagreb index with the second degree, and proposed a new index about Zagreb index, defined as the first leap Zagreb index \cite{gut17}. The first leap Zagreb index is defined as: $$ LM_{1}(G)=\sum_{v\in V(G)}d_{2}^{2}(v|G)$$

In addition to the above mentioned degree-based topological indexes, some distance-based topological indexes have also caused extensive research. In 2004, Dankelmann introduced the eccentricity sum index \cite{dan04}, defined as: $$\theta(G)=\sum_{v\in V(G)}e(v|G) $$

Sharma proposed the eccentric connectivity index  \cite{sha97}, defined as: $$\xi^{C}(G)=\sum_{v\in V(G)}d_{1}(v|G)e(v|G)$$

Recently, Naji proposed the leap eccentric connectivity index \cite{naj}, defined as:$$L\xi^{C}(G)=\sum_{v\in V(G)}d_{2}(v|G)e(v|G)$$

Regarding the study of the leap eccentric connectivity index, currently Manjunathe et al. \cite{man19} gave the values of the leap eccentric connectivity index for some graph operations contains cartesian product, composition, disjunctions, symmetric difference and corona product. The leap eccentric connectivity index of some other graphs has not been studied. This paper mainly studies the upper and lower bounds of the leap eccentric connectivity index of some graph operations with subdivided edges. Next, first introduce some basic graphs.

For a given graph $G$, the line graph is denoted by $L(G)$, and its vertex set is the edge set of $G$, where two vertices are adjacent if and only if they are adjacent edges in $G$. Cvetkovic \cite{cve80} proposed four graphs based on subdivision edges, the subdivision graph $S(G)$, the $Q(G)$ graph, the $R(G)$ graph and the total graph $T(G)$. The subdivision graph of graph $G$ is represented as $S(G)$, which represents the graph obtained by adding a vertex to each edge of $G$. The vertex set of $S(G)$ can be divided into two parts: one is the vertex set $V(G)$ in the original $G$, and the other is the vertex inserted on each edge, which is the original edge set $E(G)$ of $G$. The $Q(G)$ graph of $G$ is a graph obtained by inserting new vertices on each edge of $G$, and then connecting these new vertex pairs located on adjacent edges of $G$. The $R(G)$ graph of $G$ is defined as a graph obtained by adding new vertices corresponding to each edge to $G$ and connecting each new vertex with the two endpoints of its corresponding edge, which can also be understood as the graph obtained by replacing each edge with a triangle in $G$. The total graph $T(G)$ of $G$ is a graph whose vertex set is the union of the vertex set of $G$ and the edge set of $G$, and the two vertices of $T(G)$ are adjacent if and only if the corresponding elements in $G$ are adjacent or incident. In recent years, many researchers have paid attention to the subdivision join of several graphs. Indulal \cite{ind12} proposed two new joins, subdivision vertex join $G\dot{\vee}H $, subdivision edge join $G\underline{\vee}H $. F. Wen \cite{wen19} proposed a new graph operation based on subdivision and join on three graphs($G$ and $H_{1}$, $H_{2}$), called subdivision vertex edge join $G^{S}\bigtriangleup (H_{1}^{V}\vee H_{2}^{S})$, where $G$ and $H_{1}$, $H_{2}$ vertex sets disjoint. Barik \cite{bar17} proposed other variants of the corona graph, such as: subdivision double corona, $Q(G)$ graph double corona, $R(G)$  graph double corona, total graph double corona.

 In this paper, we first give a counterexample for that if $G$ be a graph and $S(G)$ be its  the subdivision graph, then each vertex $v\in V(G)$, $e(v|S(G))=2e(v|G)$ by Yarahmadi in \cite{yar14} in Theorem 3.1.  And we describe the upper and lower bounds of the leap eccentric connectivity index of four graphs based on subdivision edges,  the subdivision graph $S(G)$, the $Q(G)$ graph, the $R(G)$ graph and the total graph $T(G)$, and then give the expressions of the leap eccentric connectivity index of join graph based on subdivision, finally, give the bounds of the leap eccentric connectivity index of four variants of the corona graph.

\section{Leap eccentric connectivity index of four graphs based on subdivision edges}

\begin{figure}[h]
\begin{center}
  \includegraphics[width=12cm,height=3cm]{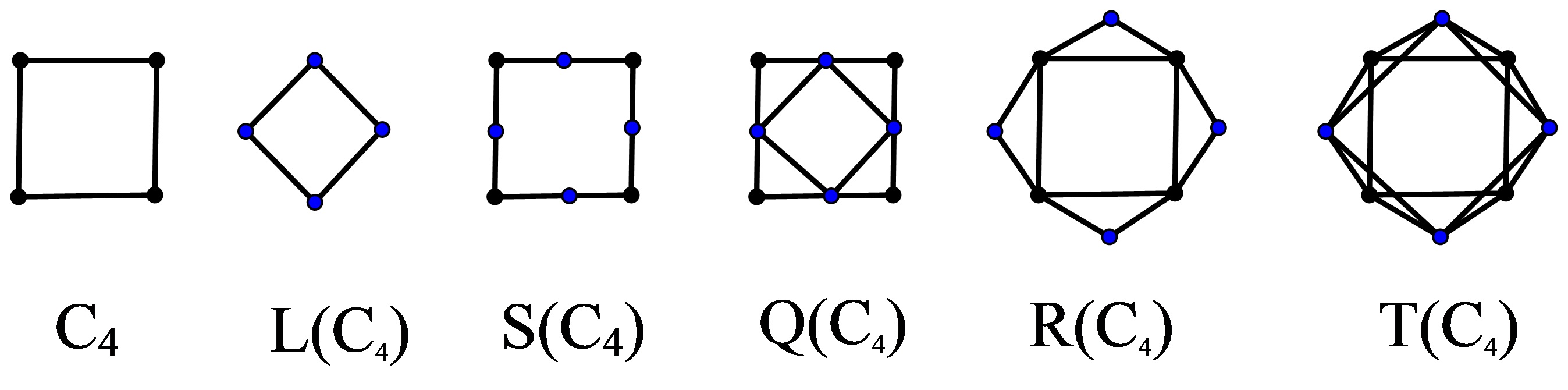}
     \end{center}
\vskip -0.5cm
\caption{{\bf Graphs $C_{4}$, $L(C_{4})$, $S(C_{4})$, $Q(C_{4})$, $R(C_{4})$, $T(C_{4})$}}\label{fig-1}
\end{figure}

This section mainly introduces the expressions of the leap eccentric connectivity index of subdivision graphs $S(G)$, the $Q(G)$ graph, the $R(G)$ graph and the total graph $T(G)$ based on other invariants. If $|V(G)|=n $, $|E(G)|=m$, we find that $S(G)$, $Q(G)$, $R(G)$, and $T(G)$ all have $n + m$ vertices. In these four graphs, we record the original $n$ vertices from graph $G$ as $v\in V(G)$, and insert the $m$ vertices corresponding to each edge of the graph $G$ are denoted as $e\in E(G)$. Taking the graph $G=C_{4}$ as an example. In Figure \ref{fig-1}, we describe the line graph $L(G)$, $S(G)$, $Q(G)$, $R(G)$, and $T(G)$, where black coloured vertices represent vertices from the vertex set in the original graph $G$, and blue coloured vertices represent newly inserted the vertex corresponding to each edge in $G$.

\begin{figure}[h]
\begin{center}
  \includegraphics[width=6cm,height=3cm]{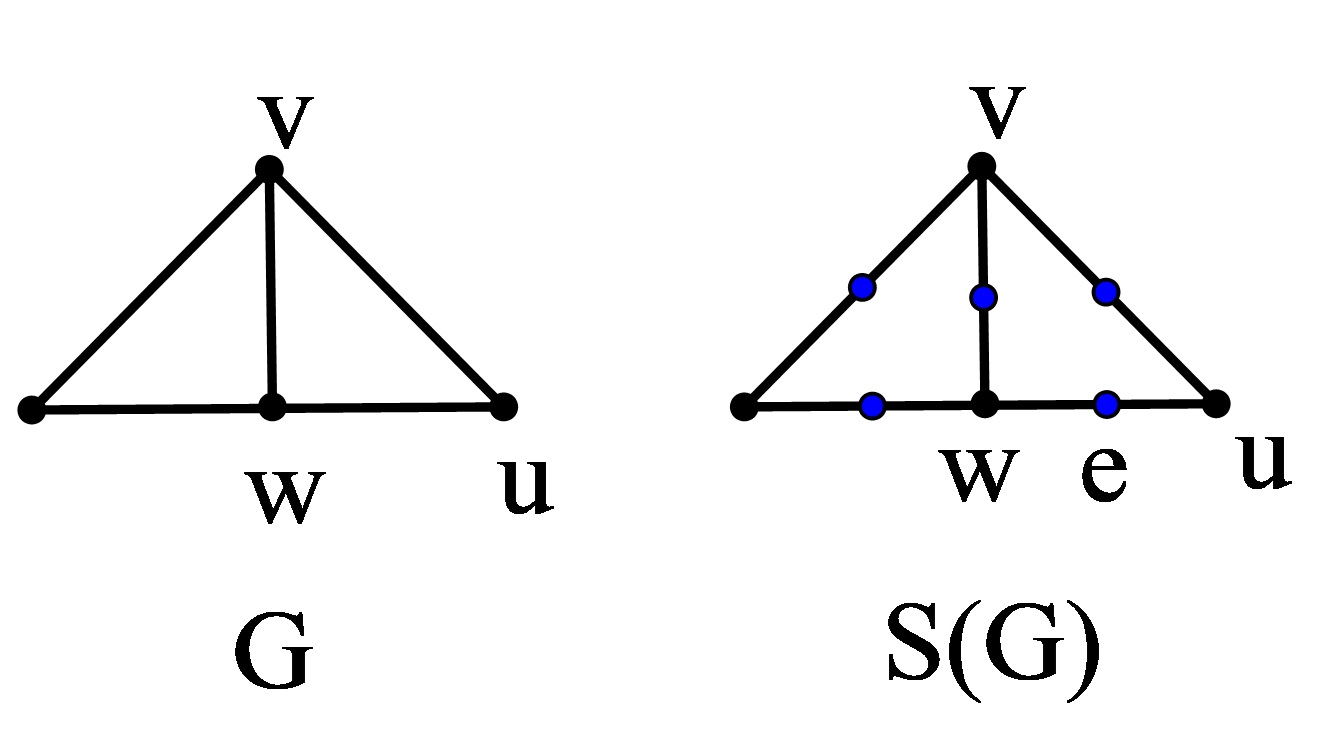}
     \end{center}
\vskip -0.5cm
\caption{{\bf A graph satisfying $e(v|S(G))=2e(v|G)+1$}}\label{fig-2}
\end{figure}

In reference \cite{yar14}, Yarahmadi gave the inequality of the eccentricity of these four graphs, but we found that $(i)$ in Theorem 3.1 for each vertex $v\in V(G)$, $e(v|S(G))=2e(v|G)$, the conclusion about the vertex eccentricity in the subdivision graph is incorrect. We can see in  Figure \ref{fig-2} that for vertex $v$, $e(v|S(G))=2e(v|G)+1$. In a graph that does not contain odd circles, $e(v|S(G))=2e(v|G)$ holds for vertex $v$.

So we correct the conclusion of Theorem 3.1 in reference \cite{yar14} to get the following lemma.

\begin{lem}\label{lem3}
Let $G$ be a connected graph, then\\
(i)For any vertex $v\in V(G)$, $2e(v|G)\leq e(v|S(G))\leq 2e(v|G)+1$;\\
(ii)For any vertex $e\in E(G)$, $2e(e|L(G))\leq e(e|S(G))\leq 2e(e|L(G))+1$.
\end{lem}

Then we introduce some conclusions related to distance in these four graphs.
\begin{lem}\label{lem1} \cite{yan07}
For any two vertices $v,v^{'}\in V(G)$, $e,e^{'}\in E(G)$, we have
$$\frac{1}{2}d_{S(G)}(v,v^{'})=d_{Q(G)}(v,v^{'})-1=d_{R(G)}(v,v^{'})=d_{T(G)}(v,v^{'})=d_{G}(v,v^{'})$$
$$\frac{1}{2}d_{S(G)}(e,e^{'})=d_{Q(G)}(e,e^{'})=d_{R(G)}(e,e^{'})-1=d_{T(G)}(e,e^{'})=d_{L(G)}(e,e^{'})$$

\end{lem}

Next, we introduce the conclusions about the second degree of vertices.
\begin{lem}\label{lem2} \cite{son16}
Let $G$ be a connected graph with $n$ vertices and $m$ edges, then
$$\sum_{v\in V(G)} d_{2}(v|G)\leq M_{1}(G)-2m $$
with equality if and only if $G$ is a $\{C_{3},C_{4}\}$-free graph.
\end{lem}

\subsection{The subdivision graph $S(G)$}
Observing graph $S(G)$ we can get the following conclusion.

\begin{lem}\label{lem4}
Let $G$ be a connected graph, then\\
(i)For any vertex $v\in V(G)$, $d_{2}(v|S(G))=d_{1}(v|G)$;\\
(ii)For any vertex $e\in E(G)$, $d_{2}(e|S(G))=d_{1}(e|L(G))$.
\end{lem}

Next we will give the bound of the leap eccentric connectivity index of the subdivision graph.

\begin{thm}\label{thm5}
Let $G$ be a connected graph, then
$$2\xi^{C}(G)+2\xi^{C}(L(G))\leq L\xi^{C}(S(G))\leq 2\xi^{C}(G)+2|E(G)|+2\xi^{C}(L(G))+2|E(L(G))| $$
\end{thm}
\zm
From Lemma \ref{lem3} and \ref{lem4}, we can get:
\begin{equation}
\begin{split}
 L\xi^{C}(S(G))&=\sum_{v\in V(G)}d_{2}(v|S(G))e(v|S(G))+\sum_{e\in E(G)}d_{2}(e|S(G))e(e|S(G))\\
 &\geq 2\sum_{v\in V(G)}d_{1}(v|G)e(v|G)+2\sum_{e\in E(G)}d_{1}(e|L(G))e(e|L(G))\\
 &=2\xi^{C}(G)+2\xi^{C}(L(G))\nonumber
 \end{split}
\end{equation}
and
\begin{equation}
\begin{split}
 L\xi^{C}(S(G))&=\sum_{v\in V(G)}d_{2}(v|S(G))e(v|S(G))+\sum_{e\in E(G)}d_{2}(e|S(G))e(e|S(G))\\
 &\leq \sum_{v\in V(G)}d_{1}(v|G)(2e(v|G)+1)+\sum_{e\in E(G)}d_{1}(e|L(G))(2e(e|L(G))+1)\\
 &=2\xi^{C}(G)+2|E(G)|+2\xi^{C}(L(G))+2|E(L(G))|\nonumber
 \end{split}
\end{equation}
The proof is completed.
\ezm

We find that when graph $G$ does not contain odd circles, $e(v|S(G))=2e(v|G)$ holds. So for a connected bipartite graph, $e(v|S(G))=2e(v|G)$  is always true, so we can get the following conclusion.

\begin{cor}\label{cor6}
Let $G$ be a connected bipartite graph, then
$$2\xi^{C}(G)+2\xi^{C}(L(G))\leq L\xi^{C}(S(G))\leq 2\xi^{C}(G)+2\xi^{C}(L(G))+2|E(L(G))|.$$
\end{cor}

When the eccentricity of all vertices in graph $G$ is reached at the pendant vertex, $e(e|S(G))=2e(e|L(G))+1$  holds for all vertices. Then we can get a connected bipartite graph where the eccentricity of all vertices is reached at the pendant vertex, the following equation holds.

\begin{cor}\label{cor7}
Let $G$ be a connected bipartite graph and the eccentricity of all vertices are reached at the pendant vertex, then
$$ L\xi^{C}(S(G))= 2\xi^{C}(G)+2\xi^{C}(L(G))+2|E(L(G))|.$$
\end{cor}

\subsection{The $Q(G)$ graph}

\begin{lem}\label{lem8} \cite{yar14}
Let $G$ be a connected graph, then\\
(i)For any vertex $v\in V(G)$, $e(v|Q(G))=e(v|G)+1$;\\
(ii)For any vertex $e\in E(G)$, $e(e|L(G))\leq e(e|Q(G))\leq e(e|L(G))+1$.
\end{lem}
Observing graph $Q(G)$ we can get the following conclusion.

\begin{lem}\label{lem9}
Let $G$ be a connected graph, then\\
(i)For any vertex $v\in V(G)$, $d_{2}(v|Q(G))\geq d_{1}(v|G)+d_{2}(v|G)$, with equality if and only if $G$ is a $\{C_{3},C_{4}\}$-free graph;\\
(ii)For any vertex $e\in E(G)$, $d_{2}(e|Q(G))\leq d_{1}(e|L(G))+d_{2}(e|L(G))$, with equality if and only if $G$ is a $\{C_{3}\}$-free graph.
\end{lem}

\begin{thm}\label{thm10}
Let $G$ be a $\{C_{3},C_{4}\}$-free connected graph, then the following inequality is established:\\
(i) $L\xi^{C}(Q(G))\geq  \xi^{C}(G)+L\xi^{C}(G)+M_{1}(G)+\xi^{C}(L(G))+L\xi^{C}(L(G))$;\\
(ii) $L\xi^{C}(Q(G))\leq \xi^{C}(G)+L\xi^{C}(G)+M_{1}(G)+\xi^{C}(L(G))+L\xi^{C}(L(G))+M_{1}(L(G))$.
\end{thm}
\zm
$G$ is a $\{C_{3},C_{4}\}$-free connected graph. According to Lemma \ref{lem9}, for any vertex $v\in V(G)$, $d_{2}(v|Q(G))= d_{1}(v|G)+d_{2}(v|G)$ holds, for any vertex $e\in E(G)$, $d_{2}(e|Q(G))=d_{1}(e|L(G))+d_{2}(e|L(G))$ holds. And according to Lemma \ref{lem8}, we can get
\begin{equation}
\begin{split}
 L\xi^{C}(Q(G))&=\sum_{v\in V(G)}d_{2}(v|Q(G))e(v|Q(G))+\sum_{e\in E(G)}d_{2}(e|Q(G))e(e|Q(G))\\
 &\geq \sum_{v\in V(G)}[d_{1}(v|G)+d_{2}(v|G)](e(v|G)+1)\\
 &+\sum_{e\in E(G)}[d_{1}(e|L(G))+d_{2}(e|L(G))]e(e|L(G))\\
 &=\xi^{C}(G)+2|E(G)|+L\xi^{C}(G)+\sum_{v\in V(G)}d_{2}(v|G)+\xi^{C}(L(G))+L\xi^{C}(L(G))\nonumber
 \end{split}
\end{equation}
In Lemma \ref{lem2}, we know that when $G$ is a $\{C_{3},C_{4}\}$-free connected graph, there is
$$\sum_{v\in V(G)} d_{2}(v|G)= M_{1}(G)-2|E(G)|. $$

Taking this equation into the above formula, we can get
\begin{equation}
\begin{split}
 L\xi^{C}(Q(G))&\geq \xi^{C}(G)+2|E(G)|+L\xi^{C}(G)+\sum_{v\in V(G)}d_{2}(v|G)+\xi^{C}(L(G))+L\xi^{C}(L(G))\\
 &=\xi^{C}(G)+2|E(G)|+L\xi^{C}(G)+M_{1}(G)-2|E(G)|+\xi^{C}(L(G))+L\xi^{C}(L(G))\\
 &=\xi^{C}(G)+L\xi^{C}(G)+M_{1}(G)+\xi^{C}(L(G))+L\xi^{C}(L(G))\nonumber
 \end{split}
\end{equation}
Inequality (i) is proved.

Similarly, we prove the second inequality (ii).
\begin{equation}
\begin{split}
 L\xi^{C}(Q(G))&=\sum_{v\in V(G)}d_{2}(v|Q(G))e(v|Q(G))+\sum_{e\in E(G)}d_{2}(e|Q(G))e(e|Q(G))\\
 &\leq \sum_{v\in V(G)}[d_{1}(v|G)+d_{2}(v|G)](e(v|G)+1)\\
 &+\sum_{e\in E(G)}[d_{1}(e|L(G))+d_{2}(e|L(G))](e(e|L(G))+1)\\
 &=\xi^{C}(G)+2|E(G)|+L\xi^{C}(G)+\sum_{v\in V(G)}d_{2}(v|G)+\xi^{C}(L(G))+L\xi^{C}(L(G))\\
 &+2|E(L(G))|+\sum_{e\in E(G)}d_{2}(e|L(G))\\
 &=\xi^{C}(G)+2|E(G)|+L\xi^{C}(G)+M_{1}(G)-2|E(G)|+\xi^{C}(L(G))+L\xi^{C}(L(G))\\
 &+2|E(L(G))|+M_{1}(L(G))-2|E(L(G))|\\
 &=\xi^{C}(G)+L\xi^{C}(G)+M_{1}(G)+\xi^{C}(L(G))+L\xi^{C}(L(G))+M_{1}(L(G))\nonumber
\end{split}
\end{equation}
The proof is completed.
\ezm

When $G$ is a $\{C_{3},C_{4}\}$-free connected bipartite graph and the eccentricity of all vertices in graph $G$ is reached at the pendant vertex, the equal sign of the second inequality in Theorem \ref{thm10} can be equal, that is, for any vertex $e\in E(G)$, $e(e|Q(G))= e(e|L(G))+1$ established. So we can get the following corollary.

\begin{cor}\label{cor11}
 Let $G$ is a $\{C_{3},C_{4}\}$-free connected bipartite graph and the eccentricity of all vertices in graph $G$ is reached at the pendant vertex, then
$$L\xi^{C}(Q(G))= \xi^{C}(G)+L\xi^{C}(G)+M_{1}(G)+\xi^{C}(L(G))+L\xi^{C}(L(G))+M_{1}(L(G))$$
\end{cor}

\subsection{The $R(G)$ graph}

\begin{lem}\label{lem12} \cite{yar14}
Let $G$ be a connected graph, then\\
(i)For any vertex $v\in V(G)$, $e(v|G)\leq e(v|R(G))\leq e(v|G)+1$;\\
(ii)For any vertex $e\in E(G)$, $e(e|R(G))= e(e|L(G))+1$.
\end{lem}
Observing graph $R(G)$ we can get the following conclusion.

\begin{lem}\label{lem13}
Let $G$ be a connected graph, then\\
(i)For any vertex $v\in V(G)$, $d_{2}(v|R(G))\geq 2d_{2}(v|G)$, with equality if and only if $G$ is a $\{C_{3},C_{4}\}$-free graph;\\
(ii)For any vertex $e\in E(G)$, $d_{2}(e|R(G))\leq  2d_{1}(e|L(G))$, with equality if and only if $G$ is a $\{C_{3}\}$-free graph.
\end{lem}

\begin{thm}\label{thm14}
Let $G$ be a $\{C_{3},C_{4}\}$-free connected graph, then \\
(i) $L\xi^{C}(R(G))\geq 2L\xi^{C}(G)+2\xi^{C}(L(G))+4|E(L(G))|$;\\
(ii) $L\xi^{C}(R(G))\leq 2L\xi^{C}(G)+2M_{1}(G)-4|E(G)|+2\xi^{C}(L(G))+4|E(L(G))|$.
\end{thm}
\zm
$G$ is a $\{C_{3},C_{4}\}$-free connected graph. According to Lemma \ref{lem13}, for any vertex $v\in V(G)$, $d_{2}(v|R(G))= 2d_{2}(v|G)$ holds, for any vertex $e\in E(G)$, $d_{2}(e|R(G))=2d_{1}(e|L(G))$ holds. And according to Lemma \ref{lem12}, we can get
\begin{equation}
\begin{split}
 L\xi^{C}(R(G))&=\sum_{v\in V(G)}d_{2}(v|R(G))e(v|R(G))+\sum_{e\in E(G)}d_{2}(e|R(G))e(e|R(G))\\
 &\geq 2\sum_{v\in V(G)}d_{2}(v|G)e(v|G)+2\sum_{e\in E(G)}d_{1}(e|L(G))(e(e|L(G))+1)\\
 &=2L\xi^{C}(G)+2\xi^{C}(L(G))+4|E(L(G))|\nonumber
 \end{split}
\end{equation}
and
\begin{equation}
\begin{split}
 L\xi^{C}(R(G))&=\sum_{v\in V(G)}d_{2}(v|R(G))e(v|R(G))+\sum_{e\in E(G)}d_{2}(e|R(G))e(e|R(G))\\
 &\leq 2\sum_{v\in V(G)}d_{2}(v|G)(e(v|G)+1)+2\sum_{e\in E(G)}d_{1}(e|L(G))(e(e|L(G))+1)\\
 &=2L\xi^{C}(G)+2\sum_{v\in V(G)}d_{2}(v|G)+2\xi^{C}(L(G))+4|E(L(G))|\nonumber
 \end{split}
\end{equation}
By Lemma \ref{lem2}, we have
\begin{equation}
\begin{split}
 L\xi^{C}(R(G))&\leq 2L\xi^{C}(G)+2(M_{1}(G)-2|E(G)|)+2\xi^{C}(L(G))+4|E(L(G))|\\
 &=2L\xi^{C}(G)+2M_{1}(G)-4|E(G)|+2\xi^{C}(L(G))+4|E(L(G))|\nonumber
 \end{split}
\end{equation}
The proof is completed.
\ezm

When $G$ is a $\{C_{3},C_{4}\}$-free connected graph with an odd circle, the equal sign on the right side of Theorem \ref{thm14} is equal, which means that for any vertex $v\in V(G)$, $e(v|R(G))=e(v|G)+1$ holds. When $G$ is a $\{C_{3},C_{4}\}$-free connected bipartite graph, then for all vertices $v\in V(G)$, $e(v|R(G))=e(v|G)$ holds. So we can get the following corollary.

\begin{cor}\label{cor15}
 Let $G$ is a $\{C_{3},C_{4}\}$-free connected bipartite graph, then
$$L\xi^{C}(R(G))= 2L\xi^{C}(G)+2\xi^{C}(L(G))+4|E(L(G))|.$$
\end{cor}

\subsection{The $T(G)$ graph}

\begin{lem}\label{lem16} \cite{yar14}
Let $G$ be a connected graph, then\\
(i)For any vertex $v\in V(G)$, $e(v|G)\leq e(v|T(G))\leq e(v|G)+1$;\\
(ii)For any vertex $e\in E(G)$, $e(e|L(G))\leq e(e|T(G))\leq e(e|L(G))+1$.
\end{lem}
Observing graph $T(G)$ we can get the following conclusion.

\begin{lem}\label{lem17}
Let $G$ be a connected graph, then\\
(i)For any vertex $v\in V(G)$, $d_{2}(v|T(G))\geq 2d_{2}(v|G)$, with equality if and only if $G$ is a $\{C_{3},C_{4}\}$-free graph;\\
(ii)For any vertex $e\in E(G)$, $d_{2}(e|T(G))\leq d_{1}(e|L(G))+d_{2}(e|L(G))$, with equality if and only if $G$ is a $\{C_{3}\}$-free graph.
\end{lem}

\begin{thm}\label{thm18}
Let $G$ be a $\{C_{3},C_{4}\}$-free connected graph, then the following inequality is established: \\
(i) $L\xi^{C}(T(G))\geq 2L\xi^{C}(G)+\xi^{C}(L(G))+L\xi^{C}(L(G))$;\\
(ii) $L\xi^{C}(T(G))\leq 2L\xi^{C}(G)+2M_{1}(G)-4|E(G)|+\xi^{C}(L(G))+L\xi^{C}(L(G))+M_{1}(L(G))$.
\end{thm}
\zm
$G$ is a $\{C_{3},C_{4}\}$-free connected graph. According to Lemma \ref{lem17}, for any vertex $v\in V(G)$, $d_{2}(v|T(G))= 2d_{2}(v|G)$ holds, for any vertex $e\in E(G)$, $d_{2}(e|T(G))=d_{1}(e|L(G))+d_{2}(e|L(G))$ holds. And according to Lemma \ref{lem16}, we can get
\begin{equation}
\begin{split}
 L\xi^{C}(T(G))&=\sum_{v\in V(G)}d_{2}(v|T(G))e(v|T(G))+\sum_{e\in E(G)}d_{2}(e|T(G))e(e|T(G))\\
 &\geq 2\sum_{v\in V(G)}d_{2}(v|G)e(v|G)+\sum_{e\in E(G)}[d_{1}(e|L(G))+d_{2}(e|L(G))]e(e|L(G))\\
 &=2L\xi^{C}(G)+\xi^{C}(L(G))+L\xi^{C}(L(G))\nonumber
 \end{split}
\end{equation}
and
\begin{equation}
\begin{split}
 L\xi^{C}(T(G))&=\sum_{v\in V(G)}d_{2}(v|T(G))e(v|T(G))+\sum_{e\in E(G)}d_{2}(e|T(G))e(e|T(G))\\
 &\leq 2\sum_{v\in V(G)}d_{2}(v|G)(e(v|G)+1)\\
 &+\sum_{e\in E(G)}[d_{1}(e|L(G))+d_{2}(e|L(G))](e(e|L(G))+1)\\
 &=2L\xi^{C}(G)+2\sum_{v\in V(G)}d_{2}(v|G)+\xi^{C}(L(G))+L\xi^{C}(L(G))+2|E(L(G))|\\
 &+\sum_{e\in E(G)}d_{2}(e|L(G))\nonumber
 \end{split}
\end{equation}
By Lemma \ref{lem2}, we have
\begin{equation}
\begin{split}
 L\xi^{C}(T(G))&\leq 2L\xi^{C}(G)+2(M_{1}(G)-2|E(G)|)+\xi^{C}(L(G))+L\xi^{C}(L(G))+2|E(L(G))|\\
 &+M_{1}(L(G))-2|E(L(G))|\\
 &=2L\xi^{C}(G)+2M_{1}(G)-4|E(G)|+\xi^{C}(L(G))++L\xi^{C}(L(G))+M_{1}(L(G))\nonumber
 \end{split}
\end{equation}
The proof is completed.
\ezm

We find that when the graph $G$ is a $\{C_{3},C_{4}\}$-free connected graph that does not contain odd cycles, $ e(v|T(G))= e(v|G)$ holds for any vertex $v$ in $G$. So for a $\{C_{3},C_{4}\}$-free connected bipartite graph, $ e(v|T(G))= e(v|G)$ is always true, so we can get the following conclusion.

\begin{cor}\label{cor19}
 Let $G$ is a $\{C_{3},C_{4}\}$-free connected bipartite graph, then\\
 (i)$L\xi^{C}(T(G))\geq  2L\xi^{C}(G)+\xi^{C}(L(G))+L\xi^{C}(L(G));$\\
 (ii)$L\xi^{C}(T(G))\leq 2L\xi^{C}(G)+\xi^{C}(L(G))+L\xi^{C}(L(G))+M_{1}(L(G)).$
\end{cor}

When graph $G$ is a $\{C_{3},C_{4}\}$-free connected graph where the eccentricity of all vertices is reached at the pendant vertex, then $e(e|T(G))=e(e|L(G))+1$ holds for all vertices $e\in E(G)$. So we can get the following equation hold.

\begin{cor}\label{cor20}
 Let $G$ is a $\{C_{3},C_{4}\}$-free connected bipartite graph and the eccentricity of all vertices is reached at the pendant vertex, then
 $$L\xi^{C}(T(G))= 2L\xi^{C}(G)+\xi^{C}(L(G))+L\xi^{C}(L(G))+M_{1}(L(G)).$$
\end{cor}

\section{The leap eccentric connectivity index of join graph based on subdivision}
In this section, we give the expressions of the leap eccentric connectivity index of several join graph based on subdivision.
\subsection{Subdivision vertex join }
\begin{defi}\label{defi3.1}
Let $G$ and $H$ be two disjoint graphs with vertices $n$ and $n_{1}$, and edges $m$ and $m_{1}$, respectively. The subdivision vertex join of graphs $G$ and $H$ is denoted as $G\dot{\vee}H$, and is a graph obtained by connecting the $i$-th$(i=1,2...,n)$ vertex in original $G$ in $S(G)$ and all vertices of $H$.
\end{defi}

\begin{lem}\label{lem3.2}
Let $G$ and $H$ be two disjoint connected graphs with vertex numbers $n\geq 3$ and $n_{1}$, and edge numbers $m$ and $m_{1}$ respectively, then for any vertex $v\in V(G\dot{\vee}H)$, we have\\
(i)If graph $G\cong S_{n}$, then
$$e(v|G\dot{\vee}H)=\left\{ \begin{array}{cl}2or3,&if v\in V(G)\\3,&if e\in E(G)\\2,&if v\in V(H)\end{array}\right. $$
(ii)If graph $G\ncong S_{n}$, then
$$e(v|G\dot{\vee}H)=\left\{ \begin{array}{cl}3,&if v\in V(G)\\3\leq e(e|G\dot{\vee}H)\leq 4 ,&if e\in E(G)\\2,&if v\in V(H)\end{array}\right. $$
(iii)
$$d_{2}(v|G\dot{\vee}H)=\left\{ \begin{array}{cl}n-1,&if v\in V(G)\\d_{1}(e|L(G))+n_{1} ,&if e\in E(G)\\n_{1}-1-d_{1}(v|H)+m,&if v\in V(H)\end{array}\right. $$
\end{lem}

\begin{thm}\label{thm3.3}
Let $G$ and $H$ be two disjoint connected graphs with vertex numbers $n\geq 3$ and $n_{1}$, and edge numbers $m$ and $m_{1}$, and $G\cong S_{n}$, then
$$L\xi^{C}(G\dot{\vee}H)=(3n-1)(n-1)+(3n+2n_{1}-3)(n+n_{1}-2)-4m_{1} $$
\end{thm}
\zm
When the graph $G\cong S_{n}=K_{1,n-1}$, we know that there is $e(v|G\dot{\vee}H)=2$ for the original center vertex in $G$, and $e(v|G\dot{\vee}H)=3$  for the other vertices in $G$, and at this time $m=n-1$, the line graph of the star graph is complete graph $K_{n-1}$. Then,
\begin{equation}
\begin{split}
 L\xi^{C}(G\dot{\vee}H)&=\sum_{v\in V(G\dot{\vee}H)}d_{2}(v|G\dot{\vee}H)e(v|G\dot{\vee}H)\\
 &=2(n-1)+3(n-1)^{2}+3\sum_{e\in E(G)}(n-2+n_{1})\\
 &+2\sum_{v\in V(H)}(n_{1}-1-d_{1}(v|H)+n-1)\\
 &=(n-1)(2+3n-3)+3(n+n_{1}-2)(n-1)+2(n+n_{1}-2)n_{1}-4m_{1}\\
 &=(3n-1)(n-1)+(3n+2n_{1}-3)(n+n_{1}-2)-4m_{1}\nonumber
 \end{split}
\end{equation}
The proof is completed.
\ezm

\begin{thm}\label{thm3.4}
Let $G$ and $H$ be two disjoint connected graphs with vertex numbers $n\geq 3$ and $n_{1}$, and edge numbers $m$ and $m_{1}$, and $G\ncong S_{n}$, then\\
(i)$L\xi^{C}(G\dot{\vee}H)\geq 3(n-1)n+6|E(L(G))|+(2n_{1}+5m-2)n_{1}-4m_{1};$\\
(ii)$L\xi^{C}(G\dot{\vee}H)\leq 3(n-1)n+8|E(L(G))|+2(n_{1}+3m-1)n_{1}-4m_{1}.$
\end{thm}
\zm
\begin{equation}
\begin{split}
 L\xi^{C}(G\dot{\vee}H)&=\sum_{v\in V(G\dot{\vee}H)}d_{2}(v|G\dot{\vee}H)e(v|G\dot{\vee}H)\\
 &=(n-1)\sum_{v\in V(G)}3+\sum_{e\in E(G)}(d_{1}(e|L(G))+n_{1})e(e|G\dot{\vee}H)\\
 &+2\sum_{v\in V(H)}(n_{1}-1-d_{1}(v|H)+m)\nonumber
 \end{split}
\end{equation}
For any vertex of $e\in E(G)$, we have $3\leq e(e|G\dot{\vee}H)\leq 4$, from which we can calculate the upper and lower bounds of $L\xi^{C}(G\dot{\vee}H)$.
\begin{equation}
\begin{split}
 L\xi^{C}(G\dot{\vee}H)&\geq 3(n-1)n+3\sum_{e\in E(G)}(d_{1}(e|L(G))+n_{1})+2(n_{1}+m-1)n_{1}-4m_{1}\\
 &=3(n-1)n+6|E(L(G))|+3n_{1}m+2(n_{1}+m-1)n_{1}-4m_{1}\\
 &=3(n-1)n+6|E(L(G))|+(2n_{1}+5m-2)n_{1}-4m_{1}\nonumber
 \end{split}
\end{equation}
Similarly,
\begin{equation}
\begin{split}
 L\xi^{C}(G\dot{\vee}H)&\leq 3(n-1)n+4\sum_{e\in E(G)}(d_{1}(e|L(G))+n_{1})+2(n_{1}+m-1)n_{1}-4m_{1}\\
 &=3(n-1)n+8|E(L(G))|+4n_{1}m+2(n_{1}+m-1)n_{1}-4m_{1}\\
 &=3(n-1)n+8|E(L(G))|+2(n_{1}+3m-1)n_{1}-4m_{1}\nonumber
 \end{split}
\end{equation}
The proof is completed.
\ezm

When the eccentricity of all vertices in graph $G$ is reached at the pendant vertex, $e(e|G\dot{\vee}H)=3$ holds for all vertices $e\in E(G)$. So we can get the following corollary.

\begin{cor}\label{cor3.5}
 Let $G$ and $H$ be two disjoint connected graphs with vertex numbers $n\geq 3$ and $n_{1}$, and edge numbers $m$ and $m_{1}$, and $G\ncong S_{n}$ where the eccentricity of all vertices is reached at the pendant vertex, then
 $$L\xi^{C}(G\dot{\vee}H)= 3(n-1)n+6|E(L(G))|+(2n_{1}+5m-2)n_{1}-4m_{1}.$$
\end{cor}

\subsection{Subdivision edge join }
\begin{defi}\label{defi3.6}
Let $G$ and $H$ be two disjoint graphs with vertices $n$ and $n_{1}$, and edges $m$ and $m_{1}$, respectively. The subdivision edge join of graphs $G$ and $H$ is denoted as $G\underline{\vee}H$, and is a graph obtained by connecting the newly inserted $j$-th$(j=1,2..,m)$ vertex corresponding to each edge of $G$ in $S(G)$ and all vertices of $H$.
\end{defi}

\begin{lem}\label{lem3.7}
Let $G$ and $H$ be two disjoint connected graphs with vertex numbers $n\geq 3$ and $n_{1}$, and edge numbers $m$ and $m_{1}$ respectively, then for any vertex $v\in V(G\underline{\vee}H)$, we have\\
(i)If graph $G\cong S_{n}$, then
$$e(v|G\underline{\vee}H)=\left\{ \begin{array}{cl}2or4,&if v\in V(G)\\3,&if e\in E(G)\\2,&if v\in V(H)\end{array}\right. $$
(ii)If graph $G\ncong S_{n}$, then
$$e(v|G\underline{\vee}H)=\left\{ \begin{array}{cl}3,&if v\in V(G),e(v|G)=1\\4,&if v\in V(G),e(v|G)>1\\3 ,&if e\in E(G)\\2,&if v\in V(H)\end{array}\right. $$
(iii)
$$d_{2}(v|G\underline{\vee}H)=\left\{ \begin{array}{cl}d_{1}(v|G)+n_{1},&if v\in V(G)\\m-1 ,&if e\in E(G)\\n_{1}-1-d_{1}(v|H)+n,&if v\in V(H)\end{array}\right. $$
\end{lem}

\begin{thm}\label{thm3.8}
Let $G$ and $H$ be two disjoint connected graphs with vertex numbers $n\geq 3$ and $n_{1}$, and edge numbers $m$ and $m_{1}$, and $G\cong S_{n}$, then
$$L\xi^{C}(G\underline{\vee}H)=(3n+4n_{1}-2)(n-1)+2(n+n_{1}-1)(n_{1}+1)-4m_{1} $$
\end{thm}
\zm
When the graph $G\cong S_{n}=K_{1,n-1}$, we know that there is $e(v|G\underline{\vee}H)=2$ for the original center vertex in $G$, and $e(v|G\underline{\vee}H)=4$  for the other vertices in $G$, and at this time $m=n-1$, the line graph of the star graph is complete graph $K_{n-1}$. Then,
\begin{equation}
\begin{split}
 L\xi^{C}(G\underline{\vee}H)&=\sum_{v\in V(G\underline{\vee}H)}d_{2}(v|G\underline{\vee}H)e(v|G\underline{\vee}H)\\
 &=2(n+n_{1}-1)+4(n_{1}+1)(n-1)+3(n-2)(n-1)\\
 &+2\sum_{v\in V(H)}(n_{1}-1-d_{1}(v|H)+n)\\
 &=2(n+n_{1}-1)+(3n+4n_{1}-2)(n-1)+2(n+n_{1}-1)n_{1}-4m_{1}\\
 &=(3n+4n_{1}-2)(n-1)+2(n+n_{1}-1)(n_{1}+1)-4m_{1}\nonumber
 \end{split}
\end{equation}
The proof is completed.
\ezm

\begin{thm}\label{thm3.9}
Let $G$ and $H$ be two disjoint connected graphs with vertex numbers $n\geq 3$ and $n_{1}$, and edge numbers $m$ and $m_{1}$, and $G\ncong S_{n}$, then
$$L\xi^{C}(G\underline{\vee}H)=m(3m+5)+2n_{1}(n_{1}+3n-1)-4m_{1}-(n+n_{1}-1)|V_{e}^{1}(G)|.$$
\end{thm}
\zm
According to Lemma \ref{lem3.7}, for any vertex $v\in V(G)$, when $e(v|G)=1$, $e(v|G\underline{\vee}H)=3$ holds, that is, when $v\in V_{e}^{1}(G)$, $e(v|G\underline{\vee}H)=3$, and $e(v|G\underline{\vee}H)=4$ holds for other original non-full vertices in $G$. Then,
\begin{equation}
\begin{split}
 L\xi^{C}(G\underline{\vee}H)&=\sum_{v\in V(G\underline{\vee}H)}d_{2}(v|G\underline{\vee}H)e(v|G\underline{\vee}H)\\
 &=3\sum_{v\in V_{e}^{1}(G)}(d_{1}(v|G)+n_{1})+4\sum_{v\notin V_{e}^{1}(G)}(d_{1}(v|G)+n_{1})+3\sum_{e\in E(G)}(m-1)\\
 &+2\sum_{v\in V(H)}(n_{1}-1-d_{1}(v|H)+n)\\
 &=4\sum_{v\in V(G)}(d_{1}(v|G)+n_{1})-\sum_{v\in V_{e}^{1}(G)}(d_{1}(v|G)+n_{1})+3m(m-1)\\
 &+2(n_{1}-1+n)n_{1}-4m_{1}\\
 &=m(3m+5)+2n_{1}(n_{1}+3n-1)-4m_{1}-(n+n_{1}-1)|V_{e}^{1}(G)|\nonumber
 \end{split}
\end{equation}
The proof is completed.
\ezm

\subsection{Subdivision vertex edge join }
\begin{defi}\label{defi3.10}
Let $G$, $H_{1}$ and $H_{2}$ be three disjoint graphs with vertices $n$, $n_{1}$ and $n_{2}$ , and edges $m$, $m_{1}$ and $m_{2}$ respectively. The subdivision vertex edge join of graphs $G$, $H_{1}$ and $H_{2}$ is denoted as $G^{S}\bigtriangleup (H_{1}^{V}\vee H_{2}^{S})$, and is the graph obtained by connecting the $i$-th$(i=1,2...,n)$ vertex in original $G$ in $S(G)$ and all vertices of $H_{1}$ and the newly inserted $j$-th$(j=1,2..,m)$ vertex corresponding to each edge of $G$ in $S(G)$ and all vertices of $H_{2}$.
\end{defi}

We can see that the number of vertices of $G^{S}\bigtriangleup (H_{1}^{V}\vee H_{2}^{S})$ is $n+m+n_{1}+n_{2}$ and the number of edges is $2m+nn_{1}+mn_{2}+m_{1}+m_{2}$.

\begin{lem}\label{lem3.11}
Let $G$, $H_{1}$ and $H_{2}$  be the disjoint three connected graphs with vertex numbers $n\geq 3$, $n_{1}$ and $n_{2}$, and edge numbers $m$, $m_{1}$ and $m_{2}$, then for any vertex $v\in V(G^{S}\bigtriangleup (H_{1}^{V}\vee H_{2}^{S}))$, we have\\
(i)If graph $G\cong S_{n}$, then
$$e(v|G^{S}\bigtriangleup (H_{1}^{V}\vee H_{2}^{S}))=\left\{ \begin{array}{cl}2or3,&if v\in V(G)\\3,&if e\in E(G)or v\in V(H_{i})\end{array}\right. $$
(ii)If graph $G\ncong S_{n}$, then
$$e(v|G^{S}\bigtriangleup (H_{1}^{V}\vee H_{2}^{S}))=3, v\in V(G^{S}\bigtriangleup (H_{1}^{V}\vee H_{2}^{S})).$$
(iii)
$$d_{2}(v|G^{S}\bigtriangleup (H_{1}^{V}\vee H_{2}^{S}))=\left\{ \begin{array}{cl}n-1+n_{2},&if v\in V(G)\\m-1+n_{1} ,&if e\in E(G)\\n_{1}-1-d_{1}(v|H_{1})+m,&if v\in V(H_{1})\\n_{2}-1-d_{1}(v|H_{2})+n,&if v\in V(H_{2})\end{array}\right. $$
\end{lem}

\begin{thm}\label{thm3.12}
Let $G$, $H_{1}$ and $H_{2}$  be the disjoint three connected graphs with vertex numbers $n\geq 3$, $n_{1}$ and $n_{2}$, and edge numbers $m$, $m_{1}$ and $m_{2}$, and $G\cong S_{n}$, then
$$L\xi^{C}(G^{S}\bigtriangleup (H_{1}^{V}\vee H_{2}^{S}))=3(n+n_{1}-2)(n+n_{1}-1)+(n+n_{2}-1)(3n+3n_{2}-1)
 -6(m_{1}+m_{2})$$
\end{thm}
\zm
When the graph $G\cong S_{n}=K_{1,n-1}$, we know that there is $e(v|G^{S}\bigtriangleup (H_{1}^{V}\vee H_{2}^{S}))=2$ for the original center vertex in $G$, and $e(v|G^{S}\bigtriangleup (H_{1}^{V}\vee H_{2}^{S}))=3$  for the other vertices in $G$, and at this time $m=n-1$, the line graph of the star graph is complete graph $K_{n-1}$. Then,
\begin{equation}
\begin{split}
 L\xi^{C}(G^{S}\bigtriangleup (H_{1}^{V}\vee H_{2}^{S}))&=\sum_{v\in V(G^{S}\bigtriangleup (H_{1}^{V}\vee H_{2}^{S}))}d_{2}(v|G^{S}\bigtriangleup (H_{1}^{V}\vee H_{2}^{S}))e(v|G^{S}\bigtriangleup (H_{1}^{V}\vee H_{2}^{S}))\\
 &=2(n-1+n_{2})+3(n-1+n_{2})(n-1)+3\sum_{e\in E(G)}(m-1+n_{1})\\
 &+3\sum_{v\in V(H_{1})}(n_{1}-1-d_{1}(v|H_{1})+m)\\
 &+3\sum_{v\in V(H_{2})}(n_{2}-1-d_{1}(v|H_{2})+n)\\
 &=3(n-1+n_{2})n-(n-1+n_{2})+3(n-2+n_{1})(n-1)\\
 &+3(n_{1}-2+n)n_{1}-6m_{1}+3(n_{2}-1+n)n_{2}-6m_{2}\\
 &=3(n+n_{1}-2)(n+n_{1}-1)+(n+n_{2}-1)(3n+3n_{2}-1)\\
 &-6(m_{1}+m_{2})\nonumber
 \end{split}
\end{equation}
The proof is completed.
\ezm

\begin{thm}\label{thm3.13}
Let $G$, $H_{1}$ and $H_{2}$  be the disjoint three connected graphs with vertex numbers $n\geq 3$, $n_{1}$ and $n_{2}$, and edge numbers $m$, $m_{1}$ and $m_{2}$, and $G\ncong S_{n}$, then
$$L\xi^{C}(G^{S}\bigtriangleup (H_{1}^{V}\vee H_{2}^{S}))=3(n+n_{2}-1)(n+n_{2})+3(m+n_{1}-1)(m+n_{1})-6(m_{1}+m_{2})$$
\end{thm}
\zm
According to Lemma \ref{lem3.11}, we have
\begin{equation}
\begin{split}
 L\xi^{C}(G^{S}\bigtriangleup (H_{1}^{V}\vee H_{2}^{S}))&=\sum_{v\in V(G^{S}\bigtriangleup (H_{1}^{V}\vee H_{2}^{S}))}d_{2}(v|G^{S}\bigtriangleup (H_{1}^{V}\vee H_{2}^{S}))e(v|G^{S}\bigtriangleup (H_{1}^{V}\vee H_{2}^{S}))\\
 &=3\sum_{v\in V(G)}(n-1+n_{2})+3\sum_{e\in E(G)}(m-1+n_{1})\\
 &+3\sum_{v\in V(H_{1})}(n_{1}-1-d_{1}(v|H_{1})+m)\\
 &+3\sum_{v\in V(H_{2})}(n_{2}-1-d_{1}(v|H_{2})+n)\\
 &=3(n-1+n_{2})n+3(m-1+n_{1})m+3(n_{1}-1+m)n_{1}\\
 &-6m_{1}+3(n_{2}-1+n)n_{2}-6m_{2}\\
 &=3(n+n_{2}-1)(n+n_{2})+3(m+n_{1}-1)(m+n_{1})-6(m_{1}+m_{2})\nonumber
 \end{split}
\end{equation}
The proof is completed.
\ezm

\section{Leap eccentric connectivity index of four double corona graphs}
In this section, we give the bounds of the leap eccentric connectivity index of four variants of the corona graph. First, an obvious conclusion is introduced.
\begin{lem}\label{lem4.1}
Let $G$ be a connected graph, then\\
(i)
$$d_{1}(v|S(G))=\left\{ \begin{array}{cl}d_{1}(v|G),&if v\in V(G)\\2,&if e\in E(G)\end{array}\right. $$
(ii)
$$d_{1}(v|Q(G))=\left\{ \begin{array}{cl}d_{1}(v|G),&if v\in V(G)\\2+d_{1}(e|L(G)),&if e\in E(G)\end{array}\right. $$
(iii)
$$d_{1}(v|R(G))=\left\{ \begin{array}{cl}2d_{1}(v|G),&if v\in V(G)\\2,&if e\in E(G)\end{array}\right. $$
(iv)
$$d_{1}(v|T(G))=\left\{ \begin{array}{cl}2d_{1}(v|G),&if v\in V(G)\\2+d_{1}(e|L(G)),&if e\in E(G)\end{array}\right. $$
\end{lem}
\subsection{Subdivision double corona }
\begin{defi}\label{defi4.2}
Let $G$ be a connected graph with $n$ vertices and $m$ edges. Let $H_{1}$ and $H_{2}$ be graphs with $n_{1}$ and $n_{2}$ vertices, respectively. The subdivision double corona graph composed of $G$, $H_{1}$ and $H_{2}$ is denoted as $G^{(S)}\circ \{H_{1},H_{2}\} $, and is a graph obtained by taking one copy of $S(G)$, $n_{1}$ copies of $H_{1}$, $m$ copies of $H_{2}$, then connecting the $i^{th}$ vertex in original $G$ in $S(G)$ to every vertex in the $i^{th}$ copy of $H_{1}$ and the newly inserted $j^{th}$-th vertex corresponding to each edge of $G$ in $S(G)$ to every vertex in the $j^{th}$ copy of $H_{2}$.
\end{defi}

We use $H_{1}^{i}$ to represent the copied $i$-th $H_{1}$, where $i=1,2,...,n$, $H_{2}^{j}$ to represent the copied $j$-th $H_{2}$, where $j=1,2,...,m$. Obviously, the number of vertices of the subdivision double corona graph $G^{(S)}\circ \{H_{1},H_{2}\} $ is $n(n_{1}+1)+m(n_{2}+1)$.

\begin{lem}\label{lem4.3}
Let $G$ be a connected graph with $n$ vertices and $m$ edges. Let $H_{1}$ and $H_{2}$ be graphs with vertices $n_{1}$ and $n_{2}$, respectively. Then\\
(i)
$$e(v|G^{(S)}\circ \{H_{1},H_{2}\})=\left\{ \begin{array}{cl}e(v|S(G))+1,&if v\in V(G)\\e(v|S(G))+1,&if e\in E(G)\\e(v_{i}|S(G))+2,&if v\in V(H_{1}^{i})\\e(e_{j}|S(G))+2,&if v\in V(H_{2}^{j})         \end{array}\right. $$
(ii)
$$d_{2}(v|G^{(S)}\circ \{H_{1},H_{2}\})=\left\{ \begin{array}{cl}(n_{2}+1)d_{1}(v|G),&if v\in V(G)\\2n_{1}+d_{1}(e|L(G)),&if e\in E(G)\\(n_{1}-1)-d_{1}(v|H_{1})+d_{1}(v_{i}|G),&if v\in V(H_{1}^{i})\\(n_{2}-1)-d_{1}(v|H_{2})+2,&if v\in V(H_{2}^{j})  \end{array}\right. $$
\end{lem}

\begin{thm}\label{thm4.4}
Let $G$ be a connected graph with $n$ vertices and $m$ edges. Let $H_{1}$ and $H_{2}$  be graphs with vertices $n_{1}$ and $ n_{2}$, and edges $m_{1}$ and $m_{2}$, respectively. Then\\
(i)
\begin{equation}
\begin{split}
 L\xi^{C}(G^{(S)}\circ \{H_{1},H_{2}\})&\geq 2(n_{1}+n_{2}+1)\xi^{C}(G)+2(n_{1}^{2}-n_{1}-2m_{1})\theta(G)+2\xi^{C}(L(G))\\
 &+2(2n_{1}+n_{2}^{2}+n_{2}-2m_{2})\theta(L(G))+2|E(L(G))|\\
 &+2m(3n_{1}+n_{2}^{2}+2n_{2}-2m_{2}+1)+2n(n_{1}^{2}-n_{1}-2m_{1})\nonumber
 \end{split}
\end{equation}
(ii)
\begin{equation}
\begin{split}
 L\xi^{C}(G^{(S)}\circ \{H_{1},H_{2}\})&\leq 2(n_{1}+n_{2}+1)\xi^{C}(G)+2(n_{1}^{2}-n_{1}-2m_{1})\theta(G)+2\xi^{C}(L(G))\\
 &+2(2n_{1}+n_{2}^{2}+n_{2}-2m_{2})\theta(L(G))+4|E(L(G))|\\
 &+m(10n_{1}+3n_{2}^{2}+7n_{2}-6m_{2}+4)+3n(n_{1}^{2}-n_{1}-2m_{1})\nonumber
 \end{split}
\end{equation}
\end{thm}
\zm
According to Lemma \ref{lem4.3}, we can get
\begin{equation}
\begin{split}
 L\xi^{C}(G^{(S)}\circ \{H_{1},H_{2}\})&=\sum_{v\in V(G^{(S)}\circ \{H_{1},H_{2}\})}d_{2}(v|G^{(S)}\circ \{H_{1},H_{2}\})e(v|G^{(S)}\circ \{H_{1},H_{2}\})\\
 &=\sum_{v\in V(G)}(n_{2}+1)d_{1}(v|G)(e(v|S(G))+1)\\
 &+\sum_{e\in E(G)}(2n_{1}+d_{1}(e|L(G)))(e(e|S(G))+1)\\
 &+\sum_{i=1}^{n}\sum_{v\in V(H_{1}^{i})}(n_{1}-1-d_{1}(v|H_{1})+d_{1}(v_{i}|G))(e(v_{i}|S(G))+2)\\
 &+\sum_{j=1}^{m}\sum_{v\in V(H_{2}^{j})}(n_{2}-1-d_{1}(v|H_{2})+2)(e(e_{j}|S(G))+2)\nonumber
 \end{split}
\end{equation}
And according to Lemma \ref{lem3}, we know that for any vertex $v\in V(G)$, $2e(v|G)\leq e(v|S(G))\leq 2e(v|G)+1$; for any vertex $e\in E(G)$, $2e(e|L(G))\leq e(e|S(G))\leq 2e(e|L(G))+1$. So we can get the upper and lower bounds of the leap eccentric connectivity index of the subdivision double corona graph $G^{(S)}\circ \{H_{1},H_{2}\}$.
\begin{equation}
\begin{split}
 L\xi^{C}(G^{(S)}\circ \{H_{1},H_{2}\})&\geq 2(n_{2}+1)\xi^{C}(G)+2(n_{2}+1)m+4n_{1}\theta(L(G))+2n_{1}m\\
 &+2\xi^{C}(L(G))+2|E(L(G))|+2n_{1}(n_{1}-1)\theta(G)+2n_{1}(n_{1}-1)n\\
 &+2n_{1}\xi^{C}(G)+4n_{1}m-4m_{1}\theta(G)-4m_{1}n+2n_{2}(n_{2}+1)\theta(L(G))\\
 &+2n_{2}(n_{2}+1)m-4m_{2}\theta(L(G))-4m_{2}m\\
 &=2(n_{1}+n_{2}+1)\xi^{C}(G)+2(n_{1}^{2}-n_{1}-2m_{1})\theta(G)+2\xi^{C}(L(G))\\
 &+2(2n_{1}+n_{2}^{2}+n_{2}-2m_{2})\theta(L(G))+2|E(L(G))|\\
 &+2m(3n_{1}+n_{2}^{2}+2n_{2}-2m_{2}+1)+2n(n_{1}^{2}-n_{1}-2m_{1})\nonumber
 \end{split}
\end{equation}
Similarly, we can get the upper bound
\begin{equation}
\begin{split}
 L\xi^{C}(G^{(S)}\circ \{H_{1},H_{2}\})&\leq 2(n_{1}+n_{2}+1)\xi^{C}(G)+2(n_{1}^{2}-n_{1}-2m_{1})\theta(G)+2\xi^{C}(L(G))\\
 &+2(2n_{1}+n_{2}^{2}+n_{2}-2m_{2})\theta(L(G))+4|E(L(G))|\\
 &+m(10n_{1}+3n_{2}^{2}+7n_{2}-6m_{2}+4)+3n(n_{1}^{2}-n_{1}-2m_{1})\nonumber
 \end{split}
\end{equation}
The proof is completed.
\ezm

For a connected bipartite graph, $e(v|S(G))=2e(v|G)$ is always true, so we can get the following conclusion.

\begin{cor}\label{cor4.5}
Let $G$ be a connected bipartite graph with $n$ vertices and $m$ edges. Let $H_{1}$ and $H_{2}$  be graphs with vertices $n_{1}$ and $ n_{2}$, and edges $m_{1}$ and $m_{2}$, respectively. Then\\
(i)
\begin{equation}
\begin{split}
 L\xi^{C}(G^{(S)}\circ \{H_{1},H_{2}\})&\geq 2(n_{1}+n_{2}+1)\xi^{C}(G)+2(n_{1}^{2}-n_{1}-2m_{1})\theta(G)+2\xi^{C}(L(G))\\
 &+2(2n_{1}+n_{2}^{2}+n_{2}-2m_{2})\theta(L(G))+2|E(L(G))|\\
 &+2m(3n_{1}+n_{2}^{2}+2n_{2}-2m_{2}+1)+2n(n_{1}^{2}-n_{1}-2m_{1})\nonumber
 \end{split}
\end{equation}
(ii)
\begin{equation}
\begin{split}
 L\xi^{C}(G^{(S)}\circ \{H_{1},H_{2}\})&\leq 2(n_{1}+n_{2}+1)\xi^{C}(G)+2(n_{1}^{2}-n_{1}-2m_{1})\theta(G)+2\xi^{C}(L(G))\\
 &+2(2n_{1}+n_{2}^{2}+n_{2}-2m_{2})\theta(L(G))+4|E(L(G))|\\
 &+m(8n_{1}+3n_{2}^{2}+5n_{2}-6m_{2}+2)+2n(n_{1}^{2}-n_{1}-2m_{1})\nonumber
 \end{split}
\end{equation}
\end{cor}
When the eccentricity of all vertices in graph $G$ is reached at the pendant vertex, $e(e|S(G))=2e(e|L(G))+1$  holds for all vertices. Then we can get the following equation for the connected bipartite graph and the eccentricity of all vertices is reached at the pendant vertex.

\begin{cor}\label{cor4.6}
Let $G$ be a connected bipartite graph and the eccentricity of all vertices are reached at the pendant vertex, where the number of vertices of $G$ is $n$ and the number of edges is $m$. Let $H_{1}$ and $H_{2}$  be graphs with vertices $n_{1}$ and $ n_{2}$, and edges $m_{1}$ and $m_{2}$, respectively. Then
\begin{equation}
\begin{split}
 L\xi^{C}(G^{(S)}\circ \{H_{1},H_{2}\})&= 2(n_{1}+n_{2}+1)\xi^{C}(G)+2(n_{1}^{2}-n_{1}-2m_{1})\theta(G)+2\xi^{C}(L(G))\\
 &+2(2n_{1}+n_{2}^{2}+n_{2}-2m_{2})\theta(L(G))+4|E(L(G))|\\
 &+m(8n_{1}+3n_{2}^{2}+5n_{2}-6m_{2}+2)+2n(n_{1}^{2}-n_{1}-2m_{1})\nonumber
 \end{split}
\end{equation}
\end{cor}

\subsection{$Q(G)$ graph double corona}
\begin{defi}\label{defi4.7}
Let $G$ be a connected graph with $n$ vertices and $m$ edges. Let $H_{1}$ and $H_{2}$ be graphs with $n_{1}$ and $n_{2}$ vertices, respectively. The $Q(G)$ graph double corona composed of $G$, $H_{1}$ and $H_{2}$ is denoted as $G^{(Q)}\circ \{H_{1},H_{2}\} $, and is a graph obtained by taking one copy of $Q(G)$, $n_{1}$ copies of $H_{1}$, $m$ copies of $H_{2}$, then connecting the $i^{th}$ vertex in original $G$ in $Q(G)$ to every vertex in the $i^{th}$ copy of $H_{1}$ and the newly inserted $j^{th}$-th vertex corresponding to each edge of $G$ in $Q(G)$ to every vertex in the $j^{th}$ copy of $H_{2}$.
\end{defi}

We use $H_{1}^{i}$ to represent the copied $i$-th $H_{1}$, where $i=1,2,...,n$, $H_{2}^{j}$ to represent the copied $j$-th $H_{2}$, where $j=1,2,...,m$. Obviously, the number of vertices of the $Q(G)$ graph double corona $G^{(Q)}\circ \{H_{1},H_{2}\} $ is $n(n_{1}+1)+m(n_{2}+1)$.

\begin{lem}\label{lem4.8}
Let $G$ be a connected graph with $n$ vertices and $m$ edges. Let $H_{1}$ and $H_{2}$ be graphs with vertices $n_{1}$ and $n_{2}$, respectively. Then\\
(i)
$$e(v|G^{(Q)}\circ \{H_{1},H_{2}\})=\left\{ \begin{array}{cl}e(v|Q(G))+1,&if v\in V(G)\\e(v|Q(G))+1,&if e\in E(G)\\e(v_{i}|Q(G))+2,&if v\in V(H_{1}^{i})\\e(e_{j}|Q(G))+2,&if v\in V(H_{2}^{j})         \end{array}\right. $$
(ii)
$$d_{2}(v|G^{(Q)}\circ \{H_{1},H_{2}\})=\left\{ \begin{array}{cl}d_{2}(v|Q(G))+n_{2}d_{1}(v|G),&if v\in V(G)\\d_{2}(e|Q(G))+2n_{1}+n_{2}d_{1}(e|L(G)),&if e\in E(G)\\(n_{1}-1)-d_{1}(v|H_{1})+d_{1}(v_{i}|G),&if v\in V(H_{1}^{i})\\(n_{2}-1)-d_{1}(v|H_{2})+2+d_{1}(e_{j}|L(G)),&if v\in V(H_{2}^{j})  \end{array}\right. $$
\end{lem}

\begin{thm}\label{thm4.9}
Let $G$ be a $\{C_{3},C_{4}\}$-free connected graph with $n$ vertices and $m$ edges. Let $H_{1}$ and $H_{2}$  be graphs with vertices $n_{1}$ and $ n_{2}$, and edges $m_{1}$ and $m_{2}$, respectively. Then\\
(i)
\begin{equation}
\begin{split}
 L\xi^{C}(G^{(Q)}\circ \{H_{1},H_{2}\})&\geq L\xi^{C}(G)+(n_{1}+n_{2}+1)\xi^{C}(G)+(n_{1}^{2}-n_{1}-2m_{1})\theta(G)\\
 &+2M_{1}(G)+L\xi^{C}(L(G))+(2n_{2}+1)\xi^{C}(L(G))\\
 &+(2n_{1}+n_{2}^{2}+n_{2}-2m_{2})\theta(L(G))+M_{1}(L(G))+6n_{2}|E(L(G))|\\
 &+2m(4n_{1}+n_{2}^{2}+3n_{2}-2m_{2})+3n(n_{1}^{2}-n_{1}-2m_{1})\nonumber
 \end{split}
\end{equation}
(ii)
\begin{equation}
\begin{split}
 L\xi^{C}(G^{(Q)}\circ \{H_{1},H_{2}\})&\leq L\xi^{C}(G)+(n_{1}+n_{2}+1)\xi^{C}(G)+(n_{1}^{2}-n_{1}-2m_{1})\theta(G)\\
 &+2M_{1}(G)+L\xi^{C}(L(G))+(2n_{2}+1)\xi^{C}(L(G))\\
 &+(2n_{1}+n_{2}^{2}+n_{2}-2m_{2})\theta(L(G))+2M_{1}(L(G))+10n_{2}|E(L(G))|\\
 &+m(10n_{1}+3n_{2}^{2}+7n_{2}-6m_{2})+3n(n_{1}^{2}-n_{1}-2m_{1})\nonumber
 \end{split}
\end{equation}
\end{thm}
\zm
When $G$ is a $\{C_{3},C_{4}\}$-free connected graph, according to Lemma \ref{lem9}, we know that for any vertex $v\in V(G)$, $d_{2}(v|Q(G))= d_{1}(v|G)+d_{2}(v|G)$, for any vertex $e\in E(G)$, $d_{2}(e|Q(G))= d_{1}(e|L(G))+d_{2}(e|L(G))$. And according to Lemma \ref{lem2}, we have when $G$ is a $\{C_{3},C_{4}\}$-free connected graph, $$\sum_{v\in V(G)}d_{2}(v|G)=M_{1}(G)-2m.$$
In Lemma \ref{lem8}, we know that for any vertex $v\in V(G)$, $e(v|Q(G))=e(v|G)+1$, for any vertex $e\in E(G)$, $e(e|L(G))\leq e(e|Q(G))\leq e(e|L(G))+1$.

According to the above lemma and Lemma \ref{lem4.8} and the definition of the leap eccentric connectivity index, we can calculate the upper and lower bounds of $L\xi^{C}(G^{(Q)}\circ \{H_{1},H_{2}\})$.

The proof is completed.
\ezm

When $G$ is a $\{C_{3},C_{4}\}$-free connected bipartite graph and the eccentricity of all vertices in graph $G$ is reached at the pendant vertex, the equal sign of the second inequality in Theorem \ref{thm4.9} can be equal, that is, for any vertex $e\in E(G)$, $e(e|Q(G))= e(e|L(G))+1$ established. So we can get the following corollary.

\begin{cor}\label{cor4.10}
 Let $G$ is a $\{C_{3},C_{4}\}$-free connected bipartite graph where the eccentricity of all vertices in graph $G$ is reached at the pendant vertex with $n$ vertices and $m$ edges. Let $H_{1}$ and $H_{2}$  be graphs with vertices $n_{1}$ and $ n_{2}$, and edges $m_{1}$ and $m_{2}$, respectively. Then
\begin{equation}
\begin{split}
 L\xi^{C}(G^{(Q)}\circ \{H_{1},H_{2}\})&= L\xi^{C}(G)+(n_{1}+n_{2}+1)\xi^{C}(G)+(n_{1}^{2}-n_{1}-2m_{1})\theta(G)\\
 &+2M_{1}(G)+L\xi^{C}(L(G))+(2n_{2}+1)\xi^{C}(L(G))\\
 &+(2n_{1}+n_{2}^{2}+n_{2}-2m_{2})\theta(L(G))+2M_{1}(L(G))+10n_{2}|E(L(G))|\\
 &+m(10n_{1}+3n_{2}^{2}+7n_{2}-6m_{2})+3n(n_{1}^{2}-n_{1}-2m_{1})\nonumber
 \end{split}
\end{equation}
\end{cor}

\subsection{$R(G)$  graph double corona }
\begin{defi}\label{defi4.11}
Let $G$ be a connected graph with $n$ vertices and $m$ edges. Let $H_{1}$ and $H_{2}$ be graphs with $n_{1}$ and $n_{2}$ vertices, respectively. The $R(G)$ graph double corona composed of $G$, $H_{1}$ and $H_{2}$ is denoted as $G^{(R)}\circ \{H_{1},H_{2}\} $, and is a graph obtained by taking one copy of $R(G)$, $n_{1}$ copies of $H_{1}$, $m$ copies of $H_{2}$, then connecting the $i^{th}$ vertex in original $G$ in $R(G)$ to every vertex in the $i^{th}$ copy of $H_{1}$ and the newly inserted $j^{th}$-th vertex corresponding to each edge of $G$ in $R(G)$ to every vertex in the $j^{th}$ copy of $H_{2}$.
\end{defi}

We use $H_{1}^{i}$ to represent the copied $i$-th $H_{1}$, where $i=1,2,...,n$, $H_{2}^{j}$ to represent the copied $j$-th $H_{2}$, where $j=1,2,...,m$. Obviously, the number of vertices of the $R(G)$ graph double corona $G^{(R)}\circ \{H_{1},H_{2}\} $ is $n(n_{1}+1)+m(n_{2}+1)$.

\begin{lem}\label{lem4.12}
Let $G$ be a connected graph with $n$ vertices and $m$ edges. Let $H_{1}$ and $H_{2}$ be graphs with vertices $n_{1}$ and $n_{2}$, respectively. Then\\
(i)
$$e(v|G^{(R)}\circ \{H_{1},H_{2}\})=\left\{ \begin{array}{cl}e(v|R(G))+1,&if v\in V(G)\\e(v|R(G))+1,&if e\in E(G)\\e(v_{i}|R(G))+2,&if v\in V(H_{1}^{i})\\e(e_{j}|R(G))+2,&if v\in V(H_{2}^{j})         \end{array}\right. $$
(ii)
$$d_{2}(v|G^{(R)}\circ \{H_{1},H_{2}\})=\left\{ \begin{array}{cl}d_{2}(v|R(G))+n_{2}d_{1}(v|G)+n_{1}d_{1}(v|G),&if v\in V(G)\\d_{2}(e|R(G))+2n_{1},&if e\in E(G)\\(n_{1}-1)-d_{1}(v|H_{1})+2d_{1}(v_{i}|G),&if v\in V(H_{1}^{i})\\(n_{2}-1)-d_{1}(v|H_{2})+2,&if v\in V(H_{2}^{j})  \end{array}\right. $$
\end{lem}

\begin{thm}\label{thm4.13}
Let $G$ be a $\{C_{3},C_{4}\}$-free connected graph with $n$ vertices and $m$ edges. Let $H_{1}$ and $H_{2}$  be graphs with vertices $n_{1}$ and $ n_{2}$, and edges $m_{1}$ and $m_{2}$, respectively. Then\\
(i)
\begin{equation}
\begin{split}
 L\xi^{C}(G^{(R)}\circ \{H_{1},H_{2}\})&\geq 2L\xi^{C}(G)+(3n_{1}+n_{2})\xi^{C}(G)+(n_{1}^{2}-n_{1}-2m_{1})\theta(G)+2M_{1}(G)\\
 &+2\xi^{C}(L(G))+(2n_{1}+n_{2}^{2}+n_{2}-2m_{2})\theta(L(G))+8|E(L(G))|\\
 &+m(14n_{1}+3n_{2}^{2}+5n_{2}-6m_{2}-4)+2n(n_{1}^{2}-n_{1}-2m_{1})\nonumber
 \end{split}
\end{equation}
(ii)
\begin{equation}
\begin{split}
 L\xi^{C}(G^{(R)}\circ \{H_{1},H_{2}\})&\leq 2L\xi^{C}(G)+(3n_{1}+n_{2})\xi^{C}(G)+(n_{1}^{2}-n_{1}-2m_{1})\theta(G)+4M_{1}(G)\\
 &+2\xi^{C}(L(G))+(2n_{1}+n_{2}^{2}+n_{2}-2m_{2})\theta(L(G))+8|E(L(G))|\\
 &+m(20n_{1}+3n_{2}^{2}+7n_{2}-6m_{2}-8)+3n(n_{1}^{2}-n_{1}-2m_{1})\nonumber
 \end{split}
\end{equation}
\end{thm}
\zm
When $G$ is a $\{C_{3},C_{4}\}$-free connected graph, according to Lemma \ref{lem13}, we know that for any vertex $v\in V(G)$, $d_{2}(v|R(G))= 2d_{2}(v|G)$, for any vertex $e\in E(G)$, $d_{2}(e|R(G))=2d_{1}(e|L(G))$. And according to Lemma \ref{lem2}, we have when $G$ is a $\{C_{3},C_{4}\}$-free connected graph, $$\sum_{v\in V(G)}d_{2}(v|G)=M_{1}(G)-2m.$$
In Lemma \ref{lem12}, we know that for any vertex $v\in V(G)$, $e(v|G)\leq e(v|R(G))\leq e(v|G)+1$, for any vertex $e\in E(G)$, $e(e|R(G))= e(e|L(G))+1$.

According to the above lemma and Lemma \ref{lem4.12} and the definition of the leap eccentric connectivity index, we can calculate the upper and lower bounds of $L\xi^{C}(G^{(R)}\circ \{H_{1},H_{2}\})$.

The proof is completed.
\ezm

When $G$ is a $\{C_{3},C_{4}\}$-free connected graph with an odd circle, the equal sign on the right side of Theorem \ref{thm4.13} is equal, which means that for any vertex $v\in V(G)$, $e(v|R(G))=e(v|G)+1$ holds. When $G$ is a $\{C_{3},C_{4}\}$-free connected bipartite graph, then for all vertices $v\in V(G)$, $e(v|R(G))=e(v|G)$ holds. So we can get the following corollary.

\begin{cor}\label{cor4.14}
 Let $G$ is a $\{C_{3},C_{4}\}$-free connected bipartite graph with $n$ vertices and $m$ edges. Let $H_{1}$ and $H_{2}$  be graphs with vertices $n_{1}$ and $ n_{2}$, and edges $m_{1}$ and $m_{2}$, respectively. Then
\begin{equation}
\begin{split}
 L\xi^{C}(G^{(R)}\circ \{H_{1},H_{2}\})&= 2L\xi^{C}(G)+(3n_{1}+n_{2})\xi^{C}(G)+(n_{1}^{2}-n_{1}-2m_{1})\theta(G)+2M_{1}(G)\\
 &+2\xi^{C}(L(G))+(2n_{1}+n_{2}^{2}+n_{2}-2m_{2})\theta(L(G))+8|E(L(G))|\\
 &+m(14n_{1}+3n_{2}^{2}+5n_{2}-6m_{2}-4)+2n(n_{1}^{2}-n_{1}-2m_{1})\nonumber
 \end{split}
\end{equation}
\end{cor}

\subsection{Total graph double corona }
\begin{defi}\label{defi4.15}
Let $G$ be a connected graph with $n$ vertices and $m$ edges. Let $H_{1}$ and $H_{2}$ be graphs with $n_{1}$ and $n_{2}$ vertices, respectively. The $T(G)$ graph double corona composed of $G$, $H_{1}$ and $H_{2}$ is denoted as $G^{(T)}\circ \{H_{1},H_{2}\} $, and is a graph obtained by taking one copy of $T(G)$, $n_{1}$ copies of $H_{1}$, $m$ copies of $H_{2}$, then connecting the $i^{th}$ vertex in original $G$ in $T(G)$ to every vertex in the $i^{th}$ copy of $H_{1}$ and the newly inserted $j^{th}$-th vertex corresponding to each edge of $G$ in $T(G)$ to every vertex in the $j^{th}$ copy of $H_{2}$.
\end{defi}

We use $H_{1}^{i}$ to represent the copied $i$-th $H_{1}$, where $i=1,2,...,n$, $H_{2}^{j}$ to represent the copied $j$-th $H_{2}$, where $j=1,2,...,m$. Obviously, the number of vertices of the $T(G)$ graph double corona $G^{(T)}\circ \{H_{1},H_{2}\} $ is $n(n_{1}+1)+m(n_{2}+1)$.

\begin{lem}\label{lem4.16}
Let $G$ be a connected graph with $n$ vertices and $m$ edges. Let $H_{1}$ and $H_{2}$ be graphs with vertices $n_{1}$ and $n_{2}$, respectively. Then\\
(i)
$$e(v|G^{(T)}\circ \{H_{1},H_{2}\})=\left\{ \begin{array}{cl}e(v|T(G))+1,&if v\in V(G)\\e(v|T(G))+1,&if e\in E(G)\\e(v_{i}|T(G))+2,&if v\in V(H_{1}^{i})\\e(e_{j}|T(G))+2,&if v\in V(H_{2}^{j})         \end{array}\right. $$
(ii)
$$d_{2}(v|G^{(T)}\circ \{H_{1},H_{2}\})=\left\{ \begin{array}{cl}d_{2}(v|T(G))+n_{2}d_{1}(v|G)+n_{1}d_{1}(v|G),&if v\in V(G)\\d_{2}(e|T(G))+2n_{1}+n_{2}d_{1}(e|L(G)),&if e\in E(G)\\(n_{1}-1)-d_{1}(v|H_{1})+2d_{1}(v_{i}|G),&if v\in V(H_{1}^{i})\\(n_{2}-1)-d_{1}(v|H_{2})+2+d_{1}(e_{j}|L(G)),&if v\in V(H_{2}^{j})  \end{array}\right. $$
\end{lem}

\begin{thm}\label{thm4.17}
Let $G$ be a $\{C_{3},C_{4}\}$-free connected graph with $n$ vertices and $m$ edges. Let $H_{1}$ and $H_{2}$  be graphs with vertices $n_{1}$ and $ n_{2}$, and edges $m_{1}$ and $m_{2}$, respectively. Then\\
(i)
\begin{equation}
\begin{split}
 L\xi^{C}(G^{(T)}\circ \{H_{1},H_{2}\})&\geq 2L\xi^{C}(G)+(3n_{1}+n_{2})\xi^{C}(G)+2M_{1}(G)+(n_{1}^{2}-n_{1}-2m_{1})\theta(G)\\
 &+L\xi^{C}(L(G))+(2n_{2}+1)\xi^{C}(L(G))+M_{1}(L(G))\\
 &+(2n_{1}+n_{2}^{2}+n_{2}-2m_{2})\theta(L(G))+6n_{2}|E(L(G))|\\
 &+2m(6n_{1}+n_{2}^{2}+2n_{2}-2m_{2}-2)+2n(n_{1}^{2}-n_{1}-2m_{1})\nonumber
 \end{split}
\end{equation}
(ii)
\begin{equation}
\begin{split}
 L\xi^{C}(G^{(T)}\circ \{H_{1},H_{2}\})&\leq 2L\xi^{C}(G)+(3n_{1}+n_{2})\xi^{C}(G)+4M_{1}(G)+(n_{1}^{2}-n_{1}-2m_{1})\theta(G)\\
 &+L\xi^{C}(L(G))+(2n_{2}+1)\xi^{C}(L(G))+2M_{1}(L(G))\\
 &+(2n_{1}+n_{2}^{2}+n_{2}-2m_{2})\theta(L(G))+10n_{2}|E(L(G))|\\
 &+m(20n_{1}+3n_{2}^{2}+7n_{2}-6m_{2}-8)+3n(n_{1}^{2}-n_{1}-2m_{1})\nonumber
 \end{split}
\end{equation}
\end{thm}
\zm
When $G$ is a $\{C_{3},C_{4}\}$-free connected graph, according to Lemma \ref{lem17}, we know that for any vertex $v\in V(G)$, $d_{2}(v|T(G))= 2d_{2}(v|G)$, for any vertex $e\in E(G)$, $d_{2}(e|T(G))=d_{1}(e|L(G))+d_{2}(e|L(G))$. And according to Lemma \ref{lem2}, we have when $G$ is a $\{C_{3},C_{4}\}$-free connected graph, $$\sum_{v\in V(G)}d_{2}(v|G)=M_{1}(G)-2m.$$
In Lemma \ref{lem16}, we know that for any vertex $v\in V(G)$, $e(v|G)\leq e(v|T(G))\leq e(v|G)+1$, for any vertex $e\in E(G)$, $e(e|L(G))\leq e(e|R(G))\leq e(e|L(G))+1$.

According to the above lemma and Lemma \ref{lem4.16} and the definition of the leap eccentric connectivity index, we can calculate the upper and lower bounds of $L\xi^{C}(G^{(T)}\circ \{H_{1},H_{2}\})$.

The proof is completed.
\ezm

We find that when the graph $G$ is a $\{C_{3},C_{4}\}$-free connected graph that does not contain odd cycles, $ e(v|T(G))= e(v|G)$ holds for any vertex $v$ in $G$. So for a $\{C_{3},C_{4}\}$-free connected bipartite graph, $ e(v|T(G))= e(v|G)$ is always true, so we can get the following conclusion.

\begin{cor}\label{cor4.18}
 Let $G$ is a $\{C_{3},C_{4}\}$-free connected bipartite graph with $n$ vertices and $m$ edges. Let $H_{1}$ and $H_{2}$  be graphs with vertices $n_{1}$ and $ n_{2}$, and edges $m_{1}$ and $m_{2}$, respectively. Then\\
 (i)
\begin{equation}
\begin{split}
 L\xi^{C}(G^{(T)}\circ \{H_{1},H_{2}\})&\geq 2L\xi^{C}(G)+(3n_{1}+n_{2})\xi^{C}(G)+2M_{1}(G)+(n_{1}^{2}-n_{1}-2m_{1})\theta(G)\\
 &+L\xi^{C}(L(G))+(2n_{2}+1)\xi^{C}(L(G))+M_{1}(L(G))\\
 &+(2n_{1}+n_{2}^{2}+n_{2}-2m_{2})\theta(L(G))+6n_{2}|E(L(G))|\\
 &+2m(6n_{1}+n_{2}^{2}+2n_{2}-2m_{2}-2)+2n(n_{1}^{2}-n_{1}-2m_{1})\nonumber
 \end{split}
\end{equation}
(ii)
\begin{equation}
\begin{split}
 L\xi^{C}(G^{(T)}\circ \{H_{1},H_{2}\})&\leq 2L\xi^{C}(G)+(3n_{1}+n_{2})\xi^{C}(G)+2M_{1}(G)+(n_{1}^{2}-n_{1}-2m_{1})\theta(G)\\
 &+L\xi^{C}(L(G))+(2n_{2}+1)\xi^{C}(L(G))+2M_{1}(L(G))\\
 &+(2n_{1}+n_{2}^{2}+n_{2}-2m_{2})\theta(L(G))+10n_{2}|E(L(G))|\\
 &+m(14n_{1}+3n_{2}^{2}+5n_{2}-6m_{2}-4)+2n(n_{1}^{2}-n_{1}-2m_{1})\nonumber
 \end{split}
\end{equation}
\end{cor}

When graph $G$ is a $\{C_{3},C_{4}\}$-free connected graph where the eccentricity of all vertices is reached at the pendant vertex, then $e(e|T(G))=e(e|L(G))+1$ holds for all vertices $e\in E(G)$. So we can get the following equation hold.

\begin{cor}\label{cor4.19}
 Let $G$ is a $\{C_{3},C_{4}\}$-free connected bipartite graph with $n$ vertices and $m$ edges, and the eccentricity of all vertices is reached at the pendant vertex. Let $H_{1}$ and $H_{2}$  be graphs with vertices $n_{1}$ and $ n_{2}$, and edges $m_{1}$ and $m_{2}$, respectively. Then\\
 \begin{equation}
\begin{split}
 L\xi^{C}(G^{(T)}\circ \{H_{1},H_{2}\})&= 2L\xi^{C}(G)+(3n_{1}+n_{2})\xi^{C}(G)+2M_{1}(G)+(n_{1}^{2}-n_{1}-2m_{1})\theta(G)\\
 &+L\xi^{C}(L(G))+(2n_{2}+1)\xi^{C}(L(G))+2M_{1}(L(G))\\
 &+(2n_{1}+n_{2}^{2}+n_{2}-2m_{2})\theta(L(G))+10n_{2}|E(L(G))|\\
 &+m(14n_{1}+3n_{2}^{2}+5n_{2}-6m_{2}-4)+2n(n_{1}^{2}-n_{1}-2m_{1})\nonumber
 \end{split}
\end{equation}
 \end{cor}

\end{document}